\documentclass[12pt]{amsart}
\usepackage{graphicx,amssymb,stmaryrd,amsfonts,epsfig,amsthm,a4,amsmath,mathrsfs}
\usepackage{vmargin,amscd,xypic,url}
\usepackage{eufrak}
\usepackage[latin1]{inputenc}
\psfigdriver{dvips}

\input xy
\xyoption{all}

\newtheorem{thm}{Theorem}[section]
\newtheorem{cor}[thm]{Corollary}
\newtheorem{lem}[thm]{Lemma}

\newtheorem{prop}[thm]{Proposition}
\theoremstyle{definition}
\newtheorem{defn}[thm]{Definition}
\newtheorem{que}[thm]{Question}
\newtheorem{rem}[thm]{Remark}
\newtheorem{exe}[thm]{Example}
\newtheorem{obs}[thm]{Observation}
\newtheorem{fact}[thm]{Fact}
\theoremstyle{remark}

\numberwithin{equation}{section}
\newcommand{\Z}{\mathbf{Z}}

\newcommand{\N}{\mathbf{N}}

\newcommand{\R}{\mathbf{R}}
\newcommand{\C}{\mathbf{C}}

\newcommand{\tpr}{\begin{tiny}\noindent Proof:}

\newcommand{\Ker}{\textnormal{Ker}}

\newcommand{\diag}{\textnormal{diag}}

\newcommand{\bpr}{\noindent \textbf{Proof}: ~}
\newcommand{\bprr}{\noindent \textbf{Proof} }
\newcommand{\epr}{~$\blacksquare$}
\begin{document}
%20E28 Maximal subgroups
%20E22 Extensions, wreath products, and other compositions
%20F05 Generators, relations, and presentations
%20B22 Multiply transitive infinite groups

\title[Finitely presented wreath products]
{Finitely presented wreath products and double coset decompositions}
\author{Yves de Cornulier}%
\date{\today}
\address{\'Ecole Polytechnique F\'ed\'erale de Lausanne (EPFL)\\
          Institut de G\'eom\'etrie, Alg\`ebre et Topologie (IGAT)\\
          1015 Lausanne, Switzerland}
\email{decornul@clipper.ens.fr}
%\thanks{}%

% ----------------------------------------------------------------
\begin{abstract}
We characterize which permutational wreath products $G\ltimes
W^{(X)}$ are finitely presented. This occurs if and only if $G$ and
$W$ are finitely presented, $G$ acts on $X$ with finitely generated
stabilizers, and with finitely many orbits on the cartesian
square~$X^2$.

On the one hand, this extends a result of G.~Baumslag about infinite
presentation of standard wreath products; on the other hand, this
provides nontrivial examples of finitely presented groups. For
instance, we obtain two quasi-isometric finitely presented groups,
one of which is torsion-free and the other has an infinite torsion
subgroup.

Motivated by the characterization above, we discuss the following
question: which finitely generated groups can have a finitely
generated subgroup with finitely many double cosets? The discussion
involves properties related to the structure of maximal subgroups,
and to the profinite topology.
\end{abstract}
\maketitle %\tableofcontents
% ----------------------------------------------------------------

\section{Introduction}

Let $G$ be a group, and $X$ a $G$-set. Let $W$ be another group.
Then $G$ acts on the direct sum $W^{(X)}$ by permutations of
factors. The {\em (permutational) wreath product} $W\wr_X G$ is
defined to be the semidirect product $W^{(X)}\rtimes G$. When the
action of $G$ on $X$ is simply transitive, it is called the {\em
standard} wreath product (this special case is sometimes called the
wreath product) and denoted by $W\wr G$.

By a result of G. Baumslag \cite{Baumslag}, a standard wreath
product $W\wr G$ with $W\neq \{1\}$ and $G$ infinite is never
finitely presented. In contrast, permutational wreath products
provide nontrivial examples:

\begin{thm}
If $W\neq \{1\}$, the wreath product $W\wr_X G$ is finitely
presented if and only if the following conditions are satisfied

(i) $W$ and $G$ are finitely presented,

(ii) $G$ acts on $X$ with finitely generated stabilizers, and

(iii) the product action of $G$ on the cartesian square $X^2$ has
finitely many orbits.\label{thm_intro_wr_fp}
\end{thm}

Note that this result extends Baumslag's result: indeed, if $G$ acts
simply transitively on $X$, then~(iii) implies that $X$ is finite.

We indicate (see Examples \ref{ThompF/FxF}, \ref{ThompF/T},
\ref{Johnson}) groups $G$ with an infinite $G$-set $X$ satisfying
the hypotheses of Theorem \ref{thm_intro_wr_fp}, which thus provides
new examples of finitely presented groups. For instance, it allows
to prove the existence of two quasi-isometric finitely presented
groups, one of which is torsion-free and the other has an infinite
torsion subgroup (see Proposition \ref{prop:Erschler_fp}).

\medskip

A general question, motivated by Theorem \ref{thm_intro_wr_fp}, is:
what are pairs $(G,X)$ satisfying the hypotheses of Theorem
\ref{thm_intro_wr_fp}? Trivial examples are pairs $(G,X)$ where $G$
is finitely presented and $X$ a finite $G$-set, thus we focus on
nontrivial cases, namely those for which $X$ is infinite.

Section \ref{Section_sg_finite_biindex_etc} is devoted to discuss
obstructions, for a given group $G$, to the existence of an infinite
$G$-set $X$ satisfying (ii) and (iii) of Theorem
\ref{thm_intro_wr_fp}, respectively satisfying (iii). It is, in the
major part, written as a survey, including many examples. For
instance, if $G$ is a finitely generated linear solvable group,
there exists no infinite $G$-set satisfying (iii) of Theorem
\ref{thm_intro_wr_fp}; while if $G$ is a free group, there exists an
infinite $G$-set satisfying (iii) of Theorem \ref{thm_intro_wr_fp},
but none can satisfy both (ii) and (iii).

\section{Finitely presented wreath products}

\subsection{Proof of Theorem \ref{thm_intro_wr_fp}}

For completeness, we first recall the following easy result.

\begin{prop}
If $X\neq\emptyset$, the wreath product $W\wr_X G$ is finitely
generated if and only if $G$ and $W$ are finitely generated, and $G$
has a finite number of orbits on $X$.\label{fgwp}
\end{prop}
\bpr If the conditions are satisfied, and if $n$ denotes the number
of $G$-orbits in $X$, then $W\wr_X G$ can be written as a quotient
of the free product $W^{\ast n}\ast G$, where $W^{\ast n}$ denotes
the free product of $n$ copies of $W$.

Conversely, suppose that $W\wr_X G$ is finitely generated. Being a
quotient of $W\wr_X G$, $G$ is also finitely generated. Since $X$ is
non-empty, $W$ embeds in $W\wr_X G$, hence is countable. If it is
not finitely generated, it can be written as the union of a strictly
increasing sequence of subgroups~$W_n$. Therefore $W\wr_X G$ is the
union of the strictly increasing sequence of subgroups $W_n\wr_X G$,
and hence is not finitely generated.\epr

\medskip

Let us now look at a presentation for the wreath product $W\wr_X G$.
For the sake of simplicity, we first suppose that $G$ acts
transitively on $X$, so that we can write $X=G/H$. It is easy to
check that a presentation for $W\wr_{G/H} G$ is given by
\begin{equation} \langle G, W,\;\;
|\;\;[H,W]\; ,\; [W,gWg^{-1}]\;\forall g\in G-H\;
\rangle\;\;\footnote{This concise notation must be understood as:
$W\wr_{G/H} G$ is the quotient of the free product $G\ast W$ by the
given relations.}.\label{pres1}\end{equation}

Using the relation $[H,W]=\{1\}$; it is immediate that, in the
family of relations $[W,gWg^{-1}]$ with $g\in G-H$, it suffices to
take into account $g$ in $G/H-\{H\}$. In fact, we can do better: we
can take $g$ in $H\backslash G/H-\{H\}$: this is obtained by
conjugating the relation $[W,gWg^{-1}]$ by an element of $H$. With
these remarks, we can prove:

\begin{thm}
Let $G$, $W$ be finitely presented groups. Let $G$ act on a set $X$,
with finitely generated stabilizers. Suppose that the product action
of $G$ on $X^2$ has a finite number of orbits. Then $W\wr_X G$ is
finitely presented.\label{wpfp}
\end{thm}
\bpr We begin by the case when $G$ is transitive on $X$, so that we
can write $X=G/H$. Since $W$ and $H$ are finitely generated,
$[H,W]=\{1\}$ in the presentation (\ref{pres1}) reduces to a finite
number of relators. The hypothesis that the product action of $G$ on
$X^2$ has a finite number of orbits reads as: $H\backslash G/H$ is
finite. Then the result follows from the remarks above: the family
of relations $[W,gWg^{-1}]$ of the presentation (\ref{pres1})
reduces to the finite family $[W,g_iWg_i^{-1}]$, where $(g_i)$ is a
finite system of representing elements of the double classes modulo
$H$ in $G$, except the class~$\{H\}$.

%%%%%%%%%%%%%%%%%%%%%%%%%%%%%%%%%%%%%%%%%%%%%%%%%%%%%%%%%%%% representants?

\medskip

We now indicate how to deal with the case when $G$ is not
necessarily transitive on $X$, which makes no essential difference.
Choosing a point in each orbit, we can write $X=\coprod_{i\in I}
G/H_i$, where $I=G\backslash X$. For all $i\in I$, consider a copy
$W_i$ of $W$. Then it is easy to check that a presentation for
$W\wr_\alpha G$ is given by the quotient of the free product of $G$
and all $W_i$ by the relations:
$$[H_i,W_i]\;\;(i\in I),\quad
[W_i,gW_ig^{-1}]\;\;(i\in I,\;g\in G-H_i),$$
$$[W_i,gW_jg^{-1}]\;\;(i,j\in I,\;i\neq j).$$

If we forget for a few seconds the two latter families of relations,
we get the generalized free product with amalgamation
$G\ast_{(H_i)}\penalty10000(H_i\times W_i)$. Given that $I$ is
finite, that $G$ and $W$ are finitely presented, and that all $H_i$
are finitely generated, this free product with amalgamation is
clearly finitely presented.

Choose $R\subset G$ such that, for every $i,j\in I$, every double
coset $H_igH_j$ is equal to $H_ig'H_j$ for some $g'\in R$. Then the
last two families of relations follow from their subfamilies when
$g$ ranges over $R$. On the other hand, the $G$-action on $X^2$
having a finite number of orbits is equivalent to saying that all
double quotients $H_i\backslash G/H_j$ are finite, so that $R$ can
be chosen finite. Thus, since $W$ is finitely generated, these
reduce to finitely many relations.\epr

\bigskip

We are now going to show that the converse of Theorem \ref{wpfp} is
true. We need some elementary preliminaries on graph products.

Let $\Gamma$ be a graph, that is, a set $\Gamma^0=I$, whose elements
are called {\em vertices}, along with a subset $\Gamma^1$ of subsets
of cardinality two of $\Gamma^0$, called {\em edges}. For each
$i\in\Gamma^0$, let $W_i$ be a group. Following \cite{Green}, the
graph product $P=(W_i)_{i\in I}^{\langle\Gamma\rangle}$ of all $W_i$
is by definition the quotient of the free product of all $W_i$ by
the relations $[W_i,W_j]=\{1\}$ if $\{i,j\}\in \Gamma^1$. Denote by
$\sigma_i$ the obvious morphism~$W_i\to P$. Observe that if $\Gamma$
is the totally disconnected graph, then $P$ is the free product of
all $W_i$, and if $\Gamma$ is the complete graph, then $\Gamma$ is
the direct sum (sometimes called the restricted direct product) of
all~$W_i$. When all $W_i$ are equal to a single group $W$; we denote
the graph product by $W^{\langle\Gamma\rangle}$.

\begin{lem}~\begin{itemize}\item[(1)]For all $i$, $\sigma_i:W_i\to P$ is injective.
\item[(2)] If $\{i,j\}\notin\Gamma^1$, the natural morphism
$\sigma_i\ast\sigma_j\to P$ is injective. \item[(3)] If
$\{i,j\}\in\Gamma^1$, the natural morphism
$\sigma_i\times\sigma_j\to P$ is injective.\item[(4)] If
$\{i,k\}\notin\Gamma^1$ and $\{j,k\}\notin\Gamma^1$, the natural
morphism $W_i\ast W_j\ast W_k\to P$, or $(W_i\times W_j)\ast W_k\to
P$ (according as whether $\{i,j\}$ belongs to $\Gamma^1$) is
injective. \label{lem:graph_ij}\end{itemize}
\end{lem}
\bpr It suffices to observe that all these morphisms are split, as
we see by taking the quotient of $P$ by the normal subgroup
generated by all $W_\ell$ for $\ell\neq i$ (resp. for $\ell\neq
i,j$) (resp. for $\ell\neq i,j,k$).\epr

\medskip

Lemma \ref{lem:graph_ij} has the following consequence. Let
$\Gamma'$ be another graph structure on the same set of vertices:
$\Gamma'^0=\Gamma^0=I$. Suppose in addition that
$\Gamma'^1\supseteq\Gamma^1$. There is a natural morphism $p$ from
$P=(W_i)_{i\in I}^{\langle\Gamma\rangle}$ to $P'=(W_i)_{i\in
I}^{\langle\Gamma'\rangle}$, which is obviously surjective. Lemma
\ref{lem:graph_ij}(2) yields:

\begin{lem}
Suppose that $W_i\neq \{1\}$ for all $i\in I$. Then the morphism $p$
is bijective if and only if
$\Gamma'^1=\Gamma^1$.\label{lem:graphprod_bij_iff}
\end{lem}
\bpr Let $\{i,j\}$ be an edge in $\Gamma^1$. Then $[W_i,W_j]=\{1\}$
in $P$. By injectivity, we get that $[W_i,W_j]=\{1\}$ in $P'$. Since
$W_i\neq\{1\}$ and $W_j\neq\{1\}$, we obtain that $W_i$ and $W_j$
cannot generate their free product in $P'$, so that, by Lemma
\ref{lem:graph_ij}(2), $\{i,j\}\in\Gamma'^1$.\epr

\medskip

Now denote by $Q$ the kernel of the natural morphism $P=(W_i)_{i\in
I}^{\langle\Gamma\rangle}\to\bigoplus_{i\in I}W_i$. We want to show
that $Q$ often contains a free non-abelian group. Assume, from now
on, that $W_i\neq \{1\}$ for all $i$. It already follows from Lemma
\ref{lem:graphprod_bij_iff} that if $\Gamma$ is not the complete
graph, then $Q\neq\{1\}$. Now denote by $\Gamma_{\textnormal{op}}$
the complement graph; namely,
$\Gamma_{\textnormal{op}}^0=\Gamma^0=I$, and, for all $i\neq j\in
I$, $\{i,j\}\in\Gamma_{\textnormal{op}}^1$ if and only if
$\{i,j\}\notin\Gamma^1$. Note that a decomposition of $\Gamma$
(resp. $\Gamma_{\textnormal{op}}$) into connected components
corresponds to a decomposition of $P$ into a free product (resp. a
direct sum).

\begin{lem}Suppose that $W_i\neq \{1\}$ for all $i$.
The following are equivalent. \begin{itemize}\item[(i)] $Q$ does not
contain any non-abelian free subgroup. \item[(ii)] All connected
components of $\Gamma_{\textnormal{op}}$ have at most 2 elements,
and whenever $\{i,j\}$ is a 2-element connected component of
$\Gamma_{\textnormal{op}}$, then $W_i$ and $W_j$ are isomorphic to
$C_2$, the cyclic group on two elements.
\end{itemize}\label{lem:F2_in_kernel}
\end{lem}
\bpr Suppose that (i) holds. Let $J$ be the union of 1-element
connected components of $\Gamma_{\textnormal{op}}$, and $K\subset
I-J$ a subset intersecting each 2-element connected component of
$\Gamma_{\textnormal{op}}$ in exactly one element. Then $Q$ can be
identified to the kernel of the natural morphism
${D_\infty}^{(K)}\to (C_2\times C_2)^{(K)}$, where $D_\infty\simeq
C_2\ast C_2$ denotes the infinite dihedral group, and thus $Q$ is
abelian (isomorphic to $\Z^{(K)}$) and cannot contains free
subgroups.

Conversely, suppose that (ii) is satisfied.

\noindent a) Suppose that there exists a connected component of
$\Gamma_{\textnormal{op}}$ with at least 2 elements, and with at
least one element $i$ such that $W_i$ is not cyclic on two elements.
Pick $j$ such that $\{i,j\}\in\Gamma_{\textnormal{op}}^1$. The
following fact is immediate.

\begin{fact}Let $G$ be a group with at least three elements. Then
it has a subgroup isomorphic to either $\Z$, $C_p$ (the cyclic group
of prime odd order $p$), $C_4$, or $C_2\times
C_2$.\label{fact:subgroup_Z_Cp_C4}
\end{fact}

Pick any nontrivial cyclic subgroup $Z_j$ in $W_j$, and any subgroup
$Z_i$ of $W_i$ as in Fact \ref{fact:subgroup_Z_Cp_C4}. By Lemma
\ref{lem:graph_ij}, there is a natural embedding of $Z_i\ast Z_j$
into $P$, which is mapped to the abelian group $Z_i\times Z_j$ in
$\bigoplus_{i\in I}W_i$. Since $Z_i\ast Z_j$ contains a non-abelian
free subgroup, so does its derived subgroup which is contained in
$Q$, so that $Q$ contains a non-abelian free subgroup.

\noindent b) Otherwise, suppose that there exists a connected
component of $\Gamma_{\textnormal{op}}$ with at least 3 elements.
Take $i,j,k\in I$, distinct, such that $\{i,k\}$ and $\{j,k\}$
belong to $\Gamma_{\textnormal{op}}^1$. We can suppose that
$W_i,W_j,W_k$ are cyclic on two elements, otherwise we can argue as
in a). By Lemma \ref{lem:graph_ij}, we get an embedding of
$(C_2\times C_2)\ast C_2$ or $C_2\ast C_2\ast C_2$ into $P$, mapping
to the abelian subgroup $C_2\times C_2\times C_2$ in
$\bigoplus_{i\in I}W_i$. As in a), since both $(C_2\times C_2)\ast
C_2$ and $C_2\ast C_2\ast C_2$ contain non-abelian free subgroups,
we obtain that $Q$ contains a non-abelian free subgroup.\epr

\medskip

When $\Gamma$ is the totally disconnected graph, Lemma
\ref{lem:F2_in_kernel} reduces as:

\begin{lem}Let $(W_i)_{i\in I}$ be a family of nontrivial groups,
and let $Q$ be the kernel of the natural morphism from the free
product of all $W_i$ to the direct sum of all $W_i$. Suppose that
$I$ has at least 2 elements, and, if all $W_i$ are cyclic on 2
elements, that $I$ has at least 3 elements. Then $Q$ contains a
non-abelian free subgroup.\epr\label{lem:glue_free_prod}
\end{lem}

\begin{lem}
Let $X$ be a set, and $\Gamma_n$ a increasing family of graph
structures on $X$: that is, $\Gamma_n^0=X$, and
$\Gamma_n^1\subseteq\Gamma_{n+1}^1$ for all $n$. Suppose that $X$
can be written as a finite disjoint union $X=\coprod_{i=1}^k X_i$
such that, for all $n$, the complement graph
$(\Gamma_n)_{\textnormal{op}}$ can be written as a disjoint union of
subgraphs $\Lambda_{n,i}$, with $\Lambda_{n,i}^0= X_i$ and
$\Lambda_{n,i}$ has constant finite degree. Then the sequence
$(\Gamma_n)$ is eventually
constant.\label{lem:seq_gros_graphes_stationne}
\end{lem}
\bpr Let $d_{n,i}$ denote the degree of $\Lambda_{n,i}$. The
sequence $(\sum_{i=1}^kd_{n,i})_{n\in\N}$ decreases, hence is
eventually constant. Thus eventually, all sequences
$(d_{n,i})_{n\in\N}$ are constant. Observe that if
$d_{n,i}=d_{n+1,i}$, then $\Lambda_{n,i}=\Lambda_{n+1,i}$.
Accordingly, the sequence $(\Gamma_{n})$ is eventually constant.\epr

\medskip

Now suppose that all $W_i$ are equal to a single group $W\neq\{1\}$,
and suppose that a group $G$ acts on $\Gamma$, i.e. acts on
$\Gamma^0=I$ preserving $\Gamma^1$. Then the semidirect product
$W^{\langle\Gamma\rangle}\rtimes G$ is well-defined.

We have to describe, given a $G$-set, what are the graph structures
preserved by $G$. Let $X$ be a set. Define an {\em edge set} on $X$
to be a subset of $X\times X$ which is symmetric and does not
intersect the diagonal; an edge set obviously defines a structure of
graph on $X$. Suppose now that $X$ is a $G$-set. Decompose $X$ into
its $G$-orbits: $X=\coprod X_i$ ($i\in I$), and choose some
base-point $x_i$ in each $X_i$ so that we can write $X_i=G/H_i$.

\begin{lem}If $E$ is a $G$-invariant edge set on $X$, and if
$(i,j)\in I^2$, define $B_{ij}=B_{ij}(E)=\{g\in G,\, (x_i, gx_j)\in
E\}$. Then the subsets $B_{ij}\subset G$ satisfy: for all $i,j\in
I$, $B_{ij}^{-1}=B_{ji}$, $H_iB_{ij}=B_{ij}$, $H_i\cap
B_{ii}=\emptyset$.

Conversely, for every family $(V_{ij})_{i,j\in I}$ of subsets of $G$
satisfying these three conditions, there exists a unique
$G$-invariant edge set $E$ such that $V_{ij}=B_{ij}(E)$ for all
$i,j\in I$, given by $(gx_i,g'x_j)\in E$ if and only if $g^{-1}g'\in
V_{ij}$.\label{G_edge_rel=dblq}\end{lem} \bpr All verifications are
straightforward.\epr

\medskip

We can now prove the converse of Theorem \ref{wpfp}. It is
essentially contained in the following slightly stronger result:

\begin{prop}
Let $G$, $W$ be groups, and $X$ a $G$-set with finitely many orbits.
Suppose that $W\neq\{1\}$, $X\neq\emptyset$, and that one of the
following conditions is satisfied.\begin{itemize}\item[(1)]The group
$G$ has infinitely many orbits on $X^2$.\item[(2)]For some $x\in X$,
the stabilizer $G_x$ is not finitely generated.\end{itemize} Then,
for every finitely presented group mapping onto $W\wr_X G$, the
kernel contains a non-abelian free subgroup. In particular, $W\wr_X
G$ is not finitely presented.\label{prop:F2_in_kernel_wreath}
\end{prop}
\bpr We keep the notation introduced above: $X=\coprod G/H_i$.

Suppose that (1) is satisfied. Then, for some $k,\ell$,
$H_k\backslash G/H_\ell$ is infinite. Define, for $i,j\in I$,
$n\in\N$, subsets $V_{ij}^n$ of $G$ as follows.

If $k\neq\ell$, take a strictly increasing sequence $(U_n)$ of
finite subsets of $H_k\backslash G/H_\ell$ whose union is all of
$H_k\backslash G/H_\ell$. Define $V_{k\ell}^n=U_n$, $V_{\ell
k}^n=U_n^{-1}$.

If $k=\ell$, take a strictly increasing sequence $(U_n)$ of finite
subsets of $H_k\backslash G/H_k-\{H_k\}$ which are symmetric under
inversion, so that the union of all $U_n$ is all of $H_k\backslash
G/H_k-\{H_k\}$. Define $V_{kk}^n=U_n$.

In both cases, for all $i,j$ such that $\{i,j\}\neq\{k,\ell\}$,
define $V_{ij}^n$ to be all of $H_i\backslash G/H_j$ if $i\neq j$,
and $H_i\backslash G/H_i-\{H_i\}$ if $i=j$.

Let $E_n$ be the $G$-invariant edge set on $X$ corresponding, by
Lemma \ref{G_edge_rel=dblq}, to the family $(V_{ij}^n)_{i,j\in I}$,
and denote by $X_n$ the corresponding graph. Observe that $(E_n)$ is
a strictly increasing sequence of $G$-invariant edge sets whose
union is the full edge set $E_\infty=X^2-\diag(X)$. Hence, the
sequence of surjective morphisms between finitely generated groups
$W^{\langle X_n\rangle}\rtimes G\to W^{\langle
X_{n+1}\rangle}\rtimes G$ converges to $W^{\langle
X_\infty\rangle}\rtimes G=W^{(X)}\rtimes G$. This already proves
that $W^{(X)}\rtimes G$ is not finitely presented: more precisely,
if a finitely presented group maps onto $W^{(X)}\rtimes G$, then the
map factors through $W^{\langle X_n\rangle}\rtimes G$ for some $n$.

Now, if the kernel of $W^{\langle X_n\rangle}\to W^{(X)}$ does not
contain a non-abelian free subgroup, then, by Lemma
\ref{lem:F2_in_kernel}, the complement graph of $X_n$ has all its
vertices of degree at most~1. Since this degree is constant on every
$G$-orbit of $X$, the hypotheses of Lemma
\ref{lem:seq_gros_graphes_stationne} are satisfied, and thus the
sequence of graphs $(X_n)$ stabilizes, a contradiction. Therefore,
for all $n$, the kernel of $W^{\langle X_n\rangle}\rtimes G\to
W^{(X)}\rtimes G$ does not contain any non-abelian free subgroup.
Since, for every finitely presented group mapping onto $W\wr_X G$,
the map must factor through $W^{\langle X_n\rangle}\rtimes G$ for
some $n$, we obtain the desired conclusion.

\medskip

Suppose that (2) is satisfied: fix $i$ such that $H_i$ is not
finitely generated. Write $H_i$ as a strictly increasing union of
subgroups $H_{i,n}$. Define $X_n$ as the disjoint union
$\coprod_{j\neq i}G/H_j\sqcup G/H_{i,n}$, and endow it with the edge
set defined as: $x\sim y$ unless $x=y$ or $x,y\in G/H_{i,n}$ and
$x\in yH_i$. Let $Q$ be the kernel of the natural map $W^{\langle
X_n\rangle}\to W^{(X)}$. It coincides with the kernel of the natural
map from the graph product $W^{\langle G/H_{i,n}\rangle}$ to
$W^{(G/H_i)}$, and hence contains the kernel of the natural map from
the free product $W^{\ast H_i/H_{i,n}}$ to $W$. Noting that
$H_{i,n}$ has infinite index in $H_i$, by Lemma
\ref{lem:glue_free_prod}, $Q$ contains a non-abelian free subgroup.
Accordingly, the kernel of $W^{\langle X_n\rangle}\rtimes G\to
W^{(X)}\rtimes G$ also contains a non-abelian free subgroup for all
$n$, and since $W^{\langle X_n\rangle}\rtimes G$ is a sequence of
finitely generated groups converging to $W\wr_X G$, we can conclude
as we did for (1): if a finitely presented group maps onto
$W^{(X)}\rtimes G$, then the map factors through $W^{\langle
X_n\rangle}\rtimes G$ for some $n$.\epr

\begin{thm}
Let $G$, $W$ be groups. Let $G$ act on a nonempty set $X$. Suppose
that $W\wr_X G$ is finitely presented. Then $G$ and $W$ are finitely
presented, and, if $W\neq \{1\}$, then the action of $G$ on $X$ has
finitely generated stabilizers, and the product action of $G$ on
$X^2$ has a finite number of orbits.\label{converse_finite_pres}
\end{thm}
\bpr By Proposition \ref{fgwp}, $G$ and $W$ are finitely generated,
and $G$ has finitely many orbits on $X$.

Now observe that $G$ is finitely presented, since it is obtained
from $W\wr_X G=W^{(X)}\rtimes G$ by killing a finite generating
subset of $W^I\subset W^{(X)}$, where $I\subset X$ is a finite set
which contains one point in each orbit.

Suppose now that $W$ is not finitely presented. Then there is a
sequence of non-injective surjective morphisms $W_{n}\to W_{n+1}$
between finitely generated groups, whose limit is $W$. Then, the
sequence of non-injective surjective morphisms between finitely
generated groups: $W_n\wr_X G\to W_{n+1}\wr_X G$ converges to
$W\wr_X G$, contradicting that $W\wr_X G$ is finitely presented.

Now Proposition \ref{prop:F2_in_kernel_wreath} allows to
conclude.\epr

\subsection{Applications}
Our main application consists in proving that the property of being
torsion-free is not weakly geometric among finitely presented
groups. The examples of \cite{Dyubina} are standard wreath products,
so are infinitely presented.

Let $F$ be the Thompson group of the dyadic interval (see Example
\ref{ThompF/FxF}), and $F_{1/2}$ the stabilizer of $1/2$. The
homogeneous space $F/F_{1/2}$ can be identified with the set $I$ of
all dyadic numbers contained in the interval $(0,1)$, and the action
of $F$ is transitive on ordered pairs $(a,b)$, $a<b$, that is,
$F_{1/2}$ has exactly three cosets in $T$.

\begin{prop}
The finitely presented groups $\Z\wr_{I}F$ and $D_{\infty}\wr_I F$
are bi-Lipschitz-equivalent. The first is torsion-free, while the
second contains an infinite subgroup of exponent
2.\label{prop:Erschler_fp}
\end{prop}
\bpr The finite presentation follows from Theorem \ref{wpfp}. The
second assertion reduces, by Proposition \ref{Erschler}, to the fact
that $\Z$ and $D_\infty$ are bi-Lipschitz-equivalent. The last
assertion is clear.\epr

\medskip

Let $S$ be any non-abelian simple, finitely presented group
(possibly finite). Let $G$ be a finitely presented group, with an
infinite index, finitely generated subgroup $H$, such that
$H\backslash G/H$ is finite, and such that the action of $G$ on
$G/H$ is faithful.

Set $\Gamma=S\wr_{G/H} G$. This group has the following properties:

\begin{prop}
1) $\Gamma$ is finitely presented.

2) Any nontrivial normal subgroup of $\Gamma$ contains
$N=S^{(G/H)}$.
\end{prop}
\bpr 1) follows from Theorem \ref{wpfp}.

2) Since the action of $G$ on $N$ is purely outer (that is, the
morphism $G\to \textnormal{Out}(N)$ is injective), every nontrivial
normal subgroup of $G$ intersects non-trivially $N$. On the other
hand, any normal subgroup $N'$ intersecting non-trivially $N$
contains it: let us recall the standard argument. For $x\in G/H$ and
$s\in S$, denote by $\delta_x(s)$ the function $X\to S$ sending $x$
to $s$ and every $y\neq x$ to $1$. If $(s_x)_{x\in G/H}$ be a
nontrivial element in $N'\cap N$, then, taking the commutator with a
suitable $\delta_x(s)$, we obtain that $N'$ contains $\delta_x(t)$
for some $x\in G/H$ and some $1\neq t\in S$. Such an element clearly
generates $N$ as a normal subgroup.\epr

\medskip

Note that the normal subgroup lattice structure of $\Gamma$ is
obtained from that of $G$ by adding a point ``at the bottom".

\medskip

An example of a direct application of Proposition
\ref{prop:F2_in_kernel_wreath} is the following well-known result,
initially proved in \cite{Shm}.

\begin{cor}
The free $d$-solvable group $R_{d,n}$ on $n$ generators ($d,n\ge 2$)
is not finitely presented.
\end{cor}
\bpr It suffices to observe that if $A$ is a finitely presented
group which maps onto $R_{d,n}$, then $A$ contains a free subgroup
of rank two. Indeed, $R_{d,n}$ maps onto $\Z\wr\Z$, while every
finitely presented group mapping onto $\Z\wr\Z$ must contain a free
subgroup by Proposition~\ref{prop:F2_in_kernel_wreath}.\epr

\section{Subgroups of finite biindex and related
properties}\label{Section_sg_finite_biindex_etc}

\subsection{Definitions and examples}

Theorems \ref{wpfp} and \ref{converse_finite_pres} raise the
following question: which finitely presented groups $G$ have an
infinite index finitely generated subgroup $H$ such that $G$ acts on
$(G/H)^2$ with a finite number of orbits? It is also natural to ask
the same question without assuming $H$ finitely generated. These
questions seem to have never been systematically investigated, but
related properties give useful information for our purposes; for
instance subgroup separability, which has been extensively studied
for other motivations, such as the generalized word problem. Hence,
the purpose of the following definitions is to present various
obstructions for a group $G$ to have an almost 2-transitive action
on an infinite set.

\begin{defn}Define a {\em pair of groups} as a pair $(G,H)$,
where $G$ is a group and $H$ a subgroup.

We say that a pair is {\em finitely presented} if $G$ is finitely
presented {\em and} $H$ is finitely generated.

We say $H$ has {\em finite biindex} in $G$ if $H\backslash G/H$ is
finite. We also say that the pair $(G,H)$ is {\em almost
2-transitive}; this is equivalent to say that $G$ has finitely many
orbits on $(G/H)^2$.

We say $H$ is {\em almost maximal} in $G$ if there are only finitely
many subgroups of $G$ containing $H$. We also say that the pair
$(G,H)$ is {\em almost primitive}.

We say that a subgroup $H$ of $G$ has {\em finite proindex} if the
profinite closure of $H$ in $G$ (that is, the intersection of all
finite index subgroups of $G$ containing $H$) has finite index in
$G$.\end{defn}

\begin{lem}
For pairs $(G,H)$, we have the implications: ($H$ has finite index)
$\Rightarrow$ ($H$ has finite biindex) $\Rightarrow$ ($H$ is almost
maximal) $\Rightarrow$ ($H$ has finite
proindex).\label{lem:implic_biindex_etc}\end{lem} \bpr The first one
is trivial. For the second one, suppose that $H$ has finite biindex
$m$. Every subgroup containing $H$ is an union of double cosets of
$H$; accordingly the number of possible subgroups is bounded by
$2^m$. For the third implication, observe that if a group has
profinite closure of infinite index, this profinite closure must be
the intersection of infinitely many finite index subgroups.\epr

\begin{rem}
None of these implications is an equivalence, even when $G$ is
finitely presented.\begin{itemize}\item For examples of infinite
index subgroups of finite biindex, see Examples \ref{ThompF/FxF},
\ref{ThompF/T}, and \ref{Johnson} below.\item If $G\neq\{1\}$ has no
proper subgroup of finite index (for instance, $G$ is infinite and
simple), then $\{1\}$ has finite proindex in $G$, but is not almost
maximal. \item Recall that, for a group $G$ and a subgroup $H$, the
pair $(G,H)$ is called a Hecke pair if, for all $g\in G$, $gHg^{-1}$
and $H$ are commensurable, i.e. they have a common finite index
subgroup; equivalently this means that the orbits of $H$ in $G/H$
are finite. On the other hand, $H$ having finite biindex means that
there are finitely many such orbits. Thus if $(G,H)$ is a Hecke pair
and $H$ has infinite index, then $H$ has infinite biindex. Now it is
known that, for any prime $p$,
$(\text{SL}_2(\Z[1/p]),\text{SL}_2(\Z))$ is a Hecke pair, and that
$\text{SL}_2(\Z)$ is a maximal subgroup of infinite index in
$\text{SL}_2(\Z[1/p])$, hence also has infinite
biindex.\end{itemize}
\end{rem}

\medskip

I only know a restricted sample of faithful almost 2-transitive
finitely presented pairs.

\begin{exe}
Let $G$ be the Thompson group $F$. This is the group of piecewise
linear increasing homeomorphisms of $[0,1]$ with singularities in
$\Z[1/2]\cap [0,1]$ and slopes powers of 2. This group is finitely
presented and torsion-free, does not contain any non-abelian free
subgroup, and has simple derived subgroup (see \cite{CFP}). The
group $F$ acts on $[0,1]\cap\Z[1/2]$, fixing $0$ and $1$, and acting
transitively on pairs $(a,b)\in\Z[1/2]$ satisfying $0<a<b<1$. The
stabilizer $F_{1/2}$ of $1/2$ is easily seen to be isomorphic to
$F\times F$. So the pair $(F,F_{1/2})$ is almost 2-transitive and
finitely presented.\label{ThompF/FxF}
\end{exe}

\begin{exe}
Let $G$ be the Thompson group $T$ (see \cite{CFP}) of the circle,
which is finitely presented and simple. This is the group of
piecewise linear oriented homeomorphisms of the circle $\R/\Z$ with
singularities in $\Z[1/2]/\Z$ and slopes powers of 2. The stabilizer
$H$ of $0=1\in\R/\Z$ is isomorphic to the Thompson group $F$ of
Example \ref{ThompF/FxF}. Then $T$ acts two-transitively on
$T/F=\Z[1/2]/\Z$. \label{ThompF/T}
\end{exe}

\begin{exe}[Houghton groups]

Fix an integer $n\ge 1$. Let $\N$ denote the non-negative integers,
and set $\Omega_n=\N\times\{1,\dots,n\}$. We think at $\Omega_n$ as
the disjoint union of $n$ copies $\N_1,\dots,\N_n$ of $\N$. Let
$G_n$ be the group of all permutations $\sigma$ of $\Omega_n$ such
that, for all $i$, $\sigma(\N_i)\Delta\N_i$ is finite, and $\sigma$
is eventually a translation on~$\N_i$.

When $n=1$, $G_1$ is the group of permutations with finite support
of $\N_1$, while $G_n$ is finitely generated if $n\ge 2$ and
finitely presented if $n\ge 3$ (see \cite{Brown}; Brown attributes
the finite presentation when $n=3$ to R. Burns and D. Solitar; the
finite generation is due to Houghton). For an explicit presentation
when $n=3$, see \cite{Johnson}.

Note that, for $n\ge 2$, the derived (resp. second derived) subgroup
of $G$ coincides with the group of permutations (resp. even
permutations) with finite support of $\Omega_n$. In particular, the
action of $G_n$ on $\Omega_n$ is $k$-transitive for all $k$.

On the other hand, as an extension of $\Z^{n-1}$ by a locally finite
group, $G_n$ is elementary amenable (but not virtually solvable).
The stabilizer $H_n$ of a point is isomorphic to $G_n$; in
particular, it is finitely generated for $n\ge 2$. \label{Johnson}
\end{exe}

\begin{exe} In \cite{RW}, a 3-manifold group $\Gamma$ together
with an infinite index surface subgroup $\Lambda$ are exhibited; it
is proved in \cite{NW} that $\Lambda\backslash\Gamma/\Lambda$ is
finite. (I do not know if the $\Gamma$-action on $\Gamma/\Lambda$ is
faithful.)
\end{exe}

\begin{exe}A refinement by D. Wise \cite{Wise} of a
construction of Rips shows that, for every finitely presented group
$Q$, there exists a finitely presented, residually finite,
torsion-free, $C'(1/6)$ small cancellation group $G$ and a
surjective map $p:G\to Q$, such that $\Ker(p)$ is a finitely
generated subgroup of $G$.

Accordingly, if $K$ is a finitely generated subgroup of finite
biindex and infinite index in $Q$, then $p^{-1}(K)$ is a finitely
generated subgroup of finite biindex and infinite index in $Q$.

Thus, starting from any of the above examples, we obtains examples
of almost 2-transitive finitely presented pairs $(G,H)$ with $G/H$
infinite and $G$ torsion-free, word hyperbolic, satisfying the
$C'(1/6)$ small cancellation property.\label{exe:rips_wise}
\end{exe}

\subsection{Related definitions}

We first introduce some obstructions to the existence of a infinite
index subgroup of finite biindex.

\begin{defn}
We say that a group $G$ has Property (PF) [respectively (MF), resp.
(BF)] if every finite proindex (resp. almost maximal, resp. finite
biindex) subgroup $H$ has finite index.
\end{defn}

We also recall that a group is (ERF) if every subgroup is closed for
the profinite topology (ERF stands for ``Extended Residually
Finite").

As a consequence of Lemma \ref{lem:implic_biindex_etc}, we have the
following implications.

$$\xymatrix{ ERF \ar@{=>}[r]
                 & PF \ar@{=>}[r] & MF
                 \ar@{=>}[r] & BF
   }$$

Note that these properties are inherited by quotients. Note also
that Properties ERF and PF are invariant by commensurability, and
that Property ERF is also inherited by subgroups. We show below
(Proposition \ref{prop:PF=MF}) that, for finitely generated groups,
Properties (PF) and (MF) are equivalent.

\begin{exe}
1) The Thompson group of Example \ref{ThompF/FxF} is 2-generated and
does not have Property (BF). In particular, non-abelian free groups
do not have Property (BF).

2) In \cite{MS}, it is proved that a finitely generated group which
is linear over a commutative ring, and not virtually solvable, has a
maximal subgroup of infinite index, thus does not have Property
(MF).

3) By a result of Olshanskii \cite{OlS_Bass-Lub_00}, any
non-elementary word hyperbolic group has an infinite quotient with
no proper subgroup of finite index. In particular, it has a maximal
subgroup of infinite index, hence does not satisfy Property (MF).

4) Hall \cite{Hall} has exhibited finitely generated 3-solvable
groups with infinite index maximal subgroups, hence without Property
(MF).

5) If $G$ is a virtually solvable group which is not virtually
polycyclic, then it is proved in \cite{Alperin} that $G$ has a
subgroup $H$ conjugate to a proper subgroup of itself. In
particular, $G$ is not ERF.

6) A virtually polycyclic group is ERF (Malcev \cite{Malcev}). It is
not known if there are other examples of finitely generated ERF
groups.

7) We prove (Proposition \ref{t:nilp_by_poly_MF}) that if a finitely
generated group $\Gamma$ is an extension with virtually polycyclic
quotient and nilpotent kernel, then $\Gamma$ has Property (MF). In
particular, this holds when $\Gamma$ is a linear virtually solvable
group.

8) The first Grigorchuk group $\Gamma$ has Property (PF) (Pervova
\cite{Per}, Grigorchuk and Wilson \cite{GW}). It is not ERF: indeed,
it has a subgroup isomorphic to a direct sum $\bigoplus_{n\ge 1}
C_{2^{n}}$, thus mapping onto the quasi-cyclic group $C_{2^\infty}$
which is not residually finite. Accordingly, $\Gamma$ has a subgroup
which is not ERF, hence $\Gamma$ is neither ERF. On the other hand,
it is an open question to find a group of subexponential growth
which does not have Property (PF); equivalently to find a group of
subexponential growth with a maximal subgroup of infinite
index.\label{exe:ERF_PF_MF_BF}\end{exe}

We now introduce similar obstructions to the existence of a
\textit{finitely generated} infinite index subgroup of finite
biindex.

\begin{defn}
We say that a group $G$ has Property (LPF) [respectively (LMF),
resp. (LBF)] if every finite proindex (resp. almost maximal, resp.
finite biindex) \textit{finitely generated} subgroup $H$ has finite
index.
\end{defn}

We also recall that a group is LERF if every finitely generated
subgroup is closed for the profinite topology. (LERF is also called
``subgroup separable"). In these four abbreviations, the additional
letter L stands for ``locally".

Again as a consequence of Lemma \ref{lem:implic_biindex_etc}, we
have the following implications.

$$\xymatrix{ ERF \ar@{=>}[d] \ar@{=>}[r]
                 & PF \ar@{=>}[d] \ar@{=>}[r] & MF
                 \ar@{=>}[d] \ar@{=>}[r] & BF \ar@{=>}[d]        \\
   LERF  \ar@{=>}[r]   & LPF \ar@{=>}[r] & LMF \ar@{=>}[r] & LBF
   }$$

Note that the properties in the second row are no longer inherited
by quotients: indeed, free groups are LERF (see the example below)
but do not have Property BF (see Example \ref{exe:ERF_PF_MF_BF}).
Note also that Properties LERF and LPF are invariant by
commensurability, and that Property LERF is also inherited by
subgroups.

In the literature, a group is defined to have the \textit{engulfing}
Property if every proper finitely generated subgroup is contained in
a proper finite index subgroup. Clearly, a group has Property (LPF)
if and only if all its finite index subgroups have the engulfing
Property.

\begin{exe}
1) A free non-abelian group is LERF \cite{HallJr}.

2) The first Grigorchuk group $\Gamma$ is LERF (Pervova \cite{Per},
Grigorchuk and Wilson \cite{GW}). On the other hand, non-residually
finite groups of subexponential growth appear in \cite{Erschler04}.

3) If $G$ is the Baumslag-Solitar group $BS(1,p)$ ($|p|\ge 2$), then
$G$ is not LERF, since its subgroup $\Z[1/p]$ is not LERF (its
quotient by a cyclic subgroup is divisible). If $G$ is a standard
wreath product $A\wr\Z$, with $A$ finitely generated abelian, then
$G$ is LERF but not ERF (Proposition \ref{wr_LERF}).

4) In \cite{NW}, an example of a free-by-cyclic 3-manifold group
which fails to satisfy Property (LPF) is given.

5) If $\Gamma\subset\text{PSL}_2(\C)$ is a lattice, then $\Gamma$
has Property (LMF). More precisely, for every finitely generated
subgroup of infinite index $\Lambda\subset\Gamma$, there exists a
strictly decreasing sequence of subgroups
$\Lambda\subset\Lambda_n\subset\Gamma$. All this follows from the
proof of \cite[Theorem 1.3]{GSS}. On the other hand, it is not known
if $\Gamma$ is always LERF, or even has (LPF).

6) Example \ref{exe:rips_wise}, and the bare existence of finitely
generated groups without Property (LBF) (Examples \ref{ThompF/FxF},
\ref{ThompF/T}, \ref{Johnson}) imply the existence of torsion-free
word hyperbolic groups without Property
(LBF).\label{exe:LERF_LPF_LMF_LBF}
\end{exe}

\subsection{Nearly maximal subgroups}

Recall \cite{Riles} that a subgroup $H$ of a group $G$ is {\it
nearly maximal} if $H$ is maximal among infinite index subgroups of
$G$. A standard verification shows that every infinite index
subgroup of a finitely generated group $G$ is contained in a nearly
maximal subgroup.

\begin{obs}
If $H$ is a nearly maximal subgroup of a group $G$, then either $H$
is closed or has finite proindex in $G$.
\end{obs}

This obvious result has the following consequence. Suppose that a
finitely generated group does not have Property (PF). Let $H$ be an
infinite index subgroup with finite proindex. Then $H$ is contained
in a nearly maximal subgroup $M$. Clearly, $M$ has also finite
proindex. So the only subgroups containing $M$ are those which
contain $\overline{M}$, and there are finitely many, so that $M$ is
almost maximal. This proves:

\begin{prop}
Let $G$ be a finitely generated group. The following are equivalent.

(i) $G$ has Property (PF).

(ii) Every nearly maximal subgroup of $G$ is profinitely closed.

(iii) $G$ has Property (MF).\label{prop:PF=MF}
\end{prop}

\begin{rem}
1) On the other hand, it is not clear whether Property (LMF) implies
Property (LPF). I actually conjecture that it is not true. A
possible counterexample could be a free product $G\ast G$, where $G$
is any nontrivial finitely generated group without any nontrivial
finite quotient, but I do not know how to prove Property (LMF) for
such a group. Note also that although $SL_n(\Z)$ $(n\ge 3)$ is known
not to have Property (LPF) \cite{SV}, whether it has Property (LMF)
is open (by \cite{MS} it does not have Property (MF)).

2) Note that the (infinitely generated) quasi-cyclic group
$C_{p^\infty}=\Z[1/p]/\Z$ is (PF) but not (MF).
\end{rem}

A consequence of Proposition \ref{prop:PF=MF} is that, for finitely
generated groups, Property (MF) is a commensurability invariant. I
do not know if this it is true for Property (LMF); however, we have:

\begin{prop}
Property (LMF) is inherited by subgroups of finite
index.\label{prop:LMF_ind_fini}
\end{prop}

\begin{lem}
Let $G$ be a group without Property (LMF). Then $G$ has a finitely
generated nearly maximal subgroup which is almost
maximal.\label{lem:LMF_nearly_max}
\end{lem}
\bpr Let $H_0$ be a finitely generated, almost maximal subgroup of
infinite index. If $H_0$ is not nearly maximal, it is properly
contained in a subgroup $K_1$ of infinite index; define $H_1$ as the
subgroup generated by $H_0$ and one element in $K_0-H_0$. Go on
defining an increasing sequence of finitely generated subgroups of
infinite index. This processus stops, since $H$ is almost maximal.
So, for some $n$, $H_n$ is nearly maximal, and, since it contains
$H_0$, it is almost maximal.\epr

\medskip

\bprr of Proposition \ref{prop:LMF_ind_fini} Note that a
non-finitely generated group necessarily has Property (LMF). Let $G$
be a finitely generated group, and $H$ a subgroup of finite index.
Suppose that $H$ does not have Property (LMF). By Lemma
\ref{lem:LMF_nearly_max}, let $M$ be a finitely generated, nearly
maximal, almost maximal subgroup of $H$. Then $M$ is contained in a
nearly maximal subgroup $M'$ of $G$. Since $M'$ has infinite index
in $G$ and $H$ has finite index in $G$, the subgroup $M'\cap H$ has
infinite index in $H$, so that $M'\cap H=M$. In particular, $M$ has
finite index in $M'$, so that $M'$ is also finitely generated.

It is clear that $M'$ is not profinitely closed in $G$: otherwise,
so would be $M=M'\cap H$, and $M$ would also be closed in $H$.\epr

\subsection{Finitely generated solvable groups}

\begin{lem}
Let $G$ be a finitely generated group which has a surjective
morphism $p$ onto an abelian group $A$, with abelian kernel $K$. Let
$H$ be a subgroup of $G$ such that $p(H)=A$. Then $H$ is closed for
the profinite topology.\label{closed_in_metab}\end{lem} \bpr The
assumption implies that $H\cap K$ is normal in $G$. Maybe replacing
$G$ by $G/(H\cap K)$, we can suppose that $H\cap K=\{1\}$, so that
$G=H\ltimes K$. To see that $H$ is closed for the profinite
topology, it clearly suffices to show that its profinite closure has
trivial intersection with $K$. Thus, let $k\in K-\{1\}$ belong to
the profinite closure of $H$.

Since $G$ is finitely generated and metabelian, it is residually
finite \cite{Hall}. So, there exists a finite index subgroup $L$ of
$K$, normal in $G$, such that $k\notin L$. Then $H\ltimes L$
contains $H$, has finite index in $G$, and does not contain $k$.
This is a contradiction.\epr

\begin{prop}
Let $G$ be a standard wreath product $A\wr\Z$, with $A\neq\{1\}$
finitely generated abelian. Then $G$ is LERF, but not
ERF.\label{wr_LERF}
\end{prop}
\bpr 1) It is not ERF because the subgroup $A^{(\N)}$ of
$A^{(\Z)}\rtimes\Z$ is not closed for the profinite topology, since
it is conjugate to a proper subgroup of itself.

2) Let $H$ be a finitely generated subgroup of $G$. Let us show that
$H$ is closed for the profinite topology.

First case: $H$ is not contained in $A^{(\Z)}$. Then the projection
of $H$ in $\Z$ is a subgroup $n\Z$ of $\Z$ ($n\ge 1$). It clearly
suffices to show that $H$ is closed in $A^{(\Z)}\rtimes n\Z$, and
this is a consequence of Lemma \ref{closed_in_metab}.

Second case: $H$ is contained in $A^{(\Z)}$. Clearly, $A^{(\Z)}$ is
closed in the profinite topology. Therefore we have to consider
$h\in A^{(\Z)}-H$ and show that $h$ is not contained in the
profinite closure of $H$. Take a finite subset $F$ of $\Z$
containing all supports of $h$ and generators of $H$. Let $n$ be
greater than the diameter of $F$. Replacing $G$ by its finite index
subgroup $A^{(\Z)}\rtimes\Z$, we can suppose that $H$ and $h$ are
contained in $A_0$, where $A^{(\Z)}=\bigoplus_{i\in\Z}A_i$. Then $h$
is a non-trivial element in the abelianization of the quotient of
$G$ by the normal subgroup generated by $H$. In particular, $h$ does
not belong to the profinite closure of $H$.\epr

% There exists a finite abelian
%group $K$ and a morphism $p_0:A_0\to K$, whose kernel contains $H$
%but not $h$. By the identification $A_i=A=A_0$, this gives rise to
%a morphism $p_i:A_i\to K$. Let $p=\bigoplus_ip_i$. Extend $p$ to
%$G$ by putting $\Z$ in the kernel. Then $p$ is a morphism of $G$
%to a finite group whose kernel contains $H$ but not $h$. Observe
%that in this special case, $h$ is not contained in the normal
%subgroup generated by $H$.

%Third case: $H$ is contained in $A^{(\Z)}$. Take a finite subset
%$F$ of $\Z$ containing all supports of generators of $H$ and $h$.
%Let $n$ be greater than the diameter of $F$. Then $H$ is contained
%in the finite index subgroup $A^{(\Z)}\rtimes n\Z$, which is
%isomorphic to $A^n\wr \Z$. This permits to reduce to the second
%case.\epr

\medskip

Example \ref{exe:LERF_LPF_LMF_LBF}(4) indicates that it is not
obvious how to generalize Proposition \ref{wr_LERF}. It would be
interesting to characterize LERF groups among finitely generated
solvable groups; even in the case of metabelian groups this is open.

\medskip

Here is now a result about Property (MF) for a class of finitely
generated solvable groups.

\begin{prop}
Let $G$ be group, which is nilpotent-by-(virtually polycyclic), i.e.
lies in a extension with nilpotent kernel and virtually polycyclic
quotient. Then $G$ has Property (MF).\label{t:nilp_by_poly_MF}
\end{prop}

Note that every finitely generated, virtually solvable group which
is linear over a field is nilpotent-by-(virtually abelian), hence
belongs to this class. In particular, this encompasses a result of
Margulis and Soifer (the easier implication in main Theorem of
\cite{MS}). The main ingredient to prove Proposition
\ref{t:nilp_by_poly_MF} is the following deep result:

\begin{thm}[Roseblade, \cite{Ros}]
Let $H$ be a virtually polycyclic group, and let $M$ be a simple $\Z
H$-module. Then $M$ is finite.\label{t:Rosenblade}
\end{thm}

\bprr of Proposition \ref{t:nilp_by_poly_MF}: let $G$ be a finitely
generated group, $N$ a nilpotent, normal subgroup, such that $G/N$
is virtually polycyclic.

Suppose by contradiction that $G$ does not have Property (MF).
Passing to a subgroup of finite index if necessary, we can suppose
that $G$ has a maximal subgroup $M$ of infinite index. We can
suppose that $M$ contains no nontrivial normal subgroup of $G$. The
centre $Z(N)$ of $N$ is normal in $G$. Since $M$ does not contain
any nontrivial normal subgroup, $M$ does not contain $Z(N)$. By
maximality, $MZ(N)=G$. Thus, since $M\cap Z(N)$ is normalized by
both $M$ and $Z(N)$, it is a normal subgroup of $G$ contained in
$M$, hence is trivial. Accordingly, $G$ is the semidirect product of
$M$ by $Z(N)$. Since $M$ is a maximal subgroup, $Z(N)$ is a simple
$M$-module, and actually a simple $M/(M\cap N)$-module since $N$
acts trivially on its centre. Since $M/(M\cap N)$ is a subgroup of
$G/N$, it is virtually polycyclic, so that by Theorem
\ref{t:Rosenblade}, $Z(N)$ is finite. Hence $M$ has finite index in
$G$, contradiction.\epr

\begin{rem}P. Hall has constructed \cite{Hall} a 3-solvable
group $G$ with a maximal, finitely generated subgroup of infinite
index. In particular, $G$ does not have Property (LMF), so that
``nilpotent-by-polycyclic" cannot be replaced by ``3-solvable" in
Proposition \ref{t:nilp_by_poly_MF}.\end{rem}

\medskip

I do not know if there exists a finitely generated solvable group
with Property (BF). However, using standard arguments, we have the
following result.

\begin{prop}
The following are equivalent.

(1) There exists a finitely generated $n$-solvable group without
Property (BF).

(2) There exists a finitely generated $n$-solvable group without
Property (LBF).

(3) There exists a finitely generated $(n-1)$-solvable group
$\Gamma$, and an infinite $\Gamma$-module $V$, such that the action
of $\Gamma$ on $V$ has finitely many orbits.\label{prop:eq_BF_solv}
\end{prop}
\bpr (2)$\Rightarrow$(1) is trivial.

(3)$\Rightarrow$(2). Observe that $\Gamma$ is a finitely generated
subgroup of finite biindex in $\Gamma\ltimes V$.

Suppose (1). Let $G$ be a finitely generated solvable group, and $M$
a subgroup of finite biindex and infinite index. Replacing $M$ by a
larger subgroup if necessary, we can suppose it nearly maximal, and
replacing $G$ by the profinite closure of $M$ if necessary, we can
suppose $M$ maximal. Moreover, taking the quotient by a normal
subgroup if necessary, we can suppose the only normal subgroup of
$G$ contained in $M$ is $\{1\}$, i.e. $G$ acts faithfully on~$G/M$.

Let $A$ be the last nontrivial term of the derived series of $G$.
Then $A$ is a normal subgroup and $A\neq \{1\}$, so that $A$ is not
contained in $M$. Accordingly, $MA=G$. Observe that $M\cap A$ is
normalized both by $M$ (since $A$ is normal) and by $A$ (since $A$
is abelian). It follows that $M\cap A$ is normal in $G$; therefore
$M\cap A=\{1\}$, and $G\simeq M\ltimes A$. Since $M$ has finite
biindex in $G$, $M$ acts with finitely many orbits on~$A$.\epr

\medskip

For $n\ge 3$, we leave as open whether the equivalent statements of
Proposition \ref{prop:eq_BF_solv} are true. For $n\le 2$, they are
false as a consequence of Proposition \ref{t:nilp_by_poly_MF}. We
record this in the following:

\begin{que}
1) Does there exist a finitely generated, solvable group without
Property (LBF)?

2) Does there exist a finitely presented solvable group without
Property (LBF)?

3) Does there exist a finitely presented solvable group without
Property (MF)?
\end{que}

The existence of a finitely generated solvable group without
Property (LBF) would permit to construct solvable finitely presented
wreath products, and would imply, arguing as in Proposition
\ref{prop:Erschler_fp}, that the class of virtually solvable groups
is not invariant under quasi-isometries within the class of finitely
presented groups.

\subsection{Amalgams and obstructions to Property ((L)BF)}\label{subs:amalgams}

The following theorem is due to M. Hall in the case of free groups,
P. Scott in the case of surface groups, and to Brunner, Burns, and
Solitar \cite{BBS} for the general case.

\begin{thm}
Let $G$ be the amalgam of two free groups over a cyclic subgroup.
Then $G$ is LERF.
\end{thm}

In contrast, Burger and Mozes \cite{BM} have constructed amalgams of
two free groups over a finite index subgroup which are finitely
presented simple groups. I do not know if these groups have Property
(LMF). These examples indicate that amalgams may have very different
behaviours, so that it seems that no general statement can be made.
The following result is a particular case of Theorem 2 in \cite{KS}.

\begin{thm}[Karrass, Solitar (1973)]
Let $G$ be a finitely generated group which splits as a nontrivial
amalgam over a finite subgroup. Then $G$ has Property (LBF).
\end{thm}

\begin{exe}
Let $G$ be an infinite group, all of whose subgroups are either
finite or of finite index. Then $G$ clearly has Property (BF). If,
moreover, $G$ has no proper subgroup of finite index, then $G$ has a
maximal subgroup which is finite; in particular $G$ does not have
Property (LMF). There exist nontrivial examples of such groups:
infinite two-generator groups all of whose nontrivial proper
subgroups are isomorphic to $\Z/p\Z$, $p$ a big prime, have been
constructed by Olshanskii (see \cite{OlS_book}). All known examples
of such groups are infinitely presented.
\end{exe}

\subsection{Fibre products}

Let $G_1,G_2,Q$ be groups and $p_i:G_i\to Q$ a surjection for
$i=1,2$. We are interested in the pair $(G,H)$, where $G=G_1\times
G_2$ and $H$ is the fibre product $G_1\times_Q G_2=\{(x,y)\in
G_1\times G_2,\; p_1(g_1)=p_2(g_2)\}$.

\begin{prop}
1) There is a natural order-preserving bijection between the set of
subgroups of $G_1\times G_2$ containing $H=G_1\times_Q G_2$ and the
set of normal subgroups of $Q$. It induces a bijection between
finite index normal subgroups of $Q$ and finite index subgroups of
$G_1\times G_2$ containing $G_1\times_Q G_2$.

2) Suppose that $G_1$ and $G_2$ are finitely generated. If $Q$ is
finitely presented, then $G_1\times_Q G_2$ is finitely generated.
Conversely, if $G_1$ and $G_2$ are in addition finitely presented
and if $G_1\times_Q G_2$ is finitely generated, then $Q$ is finitely
presented.
\end{prop}
\bpr 1) If $K$ is a subgroup of $G_1\times G_2$ containing
$H=G_1\times_Q G_2$, then set $u(K)=p_1(K\cap(G_1\times\{1\}))$.
This is a normal subgroup of $Q$, because $K\cap(G_1\times\{1\})$ is
normal in $G_1$ (identified with $G_1\times\{1\})$. Indeed, let $x$
belong to $K\cap(G_1\times\{1\})$. This means that $(x,1)\in K$. Fix
$y\in G_1$ and let us check that $yxy^{-1}\in
K\cap(G_1\times\{1\})$, i.e. $(yxy^{-1},1)\in K$. Choose $a\in G_2$
such that $p_2(a)=p_1(y)$. Then $(y,a)\in H\subseteq K$, so that
$(yxy^{-1},1)=(y,a)(x,1)(y,a)^{-1}$ also belongs to $K$.

If $N$ is a normal subgroup of $Q$, set $v(Q)=G_1\times_{Q/N} G_2$.
This is a subgroup of $G_1\times G_2$ containing $H$. We claim that
$u$ and $v$ are inverse bijections (clearly, they preserve the
order).

\begin{itemize}
    \item $K\subseteq v(u(K))$: Let $(x,y)$ belong to $K$. Write
        $p_2(y)=p_1(a)$ for some $a\in G_1$, so that $(a,y)\in H\subseteq
        K$. Then $(xa^{-1},1)=(x,y)(a,y)^{-1}\in K$. Thus
        $p_1(xa^{-1})=p_1(x)p_2(y)^{-1}\in u(K)$. This means that
        $(x,y)\in G_1\times_{Q/u(K)}G_2=v(u(K))$.
    \item $v(u(K))\subseteq K$: Let $(x,y)$ belong to $v(u(K)$.
        This means that $p_1(x)p_2(y)^{-1}\in u(K)$, i.e.
        $p_1(x)p_2(y)^{-1}=p_1(a)$ for some $a\in G_1$ such that
        $(a,1)\in K$. Therefore $(xa^{-1},y)\in H\subseteq K$, so that
        $(x,y)=(xa^{-1},y)(a,1)\in K$.
    \item $N\subseteq u(v(N))$: Let $\alpha$ belong to $N$. Write
    $\alpha=p_1(x)$ for some $x\in X$. Then $(x,1)\in
    G_1\times_{Q/N}G_1=v(N)$, so that $\alpha\in u(v(N))$.
    \item $u(v(N))\subseteq N$: Let $\alpha$ belong to $u(v(N))$.
        This means that $\alpha=p_1(x)$, for some $x\in G_1$ such
        that $(x,1)\in G_1\times_{Q/N}G_2$, so that $p_1(x)=\alpha\in N$.
\end{itemize}
 2) Suppose that $G_1$ and $G_2$ are finitely generated, and $Q$
 is finitely presented. For $i=1,2$, write $N_i=\text{Ker}(p_i)$.
 Since $Q$ is finitely presented and $G_i$ finitely generated,
 $N_i$ is generated as a normal subgroup in $G_i$ by a finite subset $R_i$.
 Besides, take a finite subset $S$ of $H=G_1\times_Q G_2$ such
 that $p_i(S)$ generates $Q$ for~$i=1,2$. Then
 $(R_1\times\{1\})\cup(\{1\}\times R_2)\cup S$ is a finite generating
 subset for a subgroup $M$ of $H$. We claim that $M=H$.
 Let $(x,y)$ belong to $H$. The hypothesis on $S$ implies that
 there exists $z\in G_2$ and $a\in N_1$ such that $(ax,z)\in M$.
 Since $(ax,z)\in M\subseteq H$, $p_1(x)=p_2(z)$, so that
 $z^{-1}y\in N_2$. Hence $(x,y)=(a,1)^{-1}(ax,z)(1,z^{-1}y)$.
We claim that $(a,1)\in M$. Indeed, $a\in N_1$, and, using that
$p_1(S)$ generates $G_1$ and $R_1$ generates $N_1$ as a normal
subgroup, we obtain that $(a,1)\in M$. Similarly, $(1,z^{-1}y)\in
M$, and therefore $(x,y)\in M$.

Conversely, suppose that $G_1$ and $G_2$ are finitely presented and
suppose that $H$ is finitely generated. There exists a finitely
presented group $Q'$, a surjective map $q:Q'\to Q$, surjective maps
$q_i:G_i\to Q'$, $i=1,2$, such that $p_i=q\circ q_i$ for $i=1,2$. If
$Q$ is not finitely presented, then the kernel of $q$ can be written
as a union of an increasing sequence of subgroups $M_n$, normal in
$Q'$. By (1), the normal subgroups intermediate between
$\text{Ker}(q_1)$ and $\text{Ker}(p_1)$ correspond bijectively with
the subgroups intermediate between $G_1\times_{Q'} G_2$ and
$G_1\times_Q G_2$. Accordingly, the latter is not finitely
generated.\epr

\begin{cor}
Let $G_1$, $G_2$ be finitely generated groups. The subgroup
$G_1\times_Q G_2$ has finite proindex (resp. is almost maximal) in
$G_1\times G_2$ if and only if $Q$ has a minimal finite index
subgroup (resp. $Q$ has a finite number of normal subgroups).\epr
\end{cor}

\begin{rem}
The question of finite presentability of a fibre product
$G_1\times_Q G_2$ is not trivial at all. It is easy to show that
$G_1$ and $G_2$, and $Q$ must necessarily be finitely presented, but
the converse is not true. For instance, take the Baumslag-Solitar
group $BS(1,p)=\Z\ltimes_p\Z[1/p]$ with $p\ge 2$, which has
presentation $\langle x,y\;|\;^xy=y^p\rangle$. There are two
morphisms $p_+$, $p_-$ of this group onto $\Z$. This gives, up to
isomorphism, two possible fiber products over $\Z$, which we denote
by $BS(1,p)\times_{\Z++}BS(1,p)$ and $BS(1,p)\times_{\Z+-}BS(1,p)$.
Then the former is finitely presented, while the second is not. The
first has presentation $\langle
x,y,z\;|\;[y,z]=1,\;^xy=y^p,\;^xz=z^p\rangle$, while the second has
a finitely generated central extension by $\Z[1/p]$, given by the
semidirect product of the diagonal subgroup of
$\textnormal{SL}_2(\Z[1/p])$ by the Heisenberg group $H_3(\Z[1/p])$,
hence is not finitely presented. For more about the finite
presentation of fibre products, see \cite{Gru,BBHM,BridGrun}.
\end{rem}

\begin{prop}
The subgroup $H=G_1\times_Q G_2$ is has finite biindex in $G_1\times
G_2$ if and only if $Q$ has a finite number of conjugacy classes.
\end{prop}
\bpr It suffices to check that every double coset of $H$ contains an
element of $G_1\times\{1\}$, and that two elements $(x,1)$ and
$(y,1)$ of $G_1\times\{1\}$ are in the same double coset if and only
if the images of $x$ and $y$ in $Q$ are conjugate.\epr

\begin{rem}
Examples of infinite, finitely generated groups with finitely many
conjugacy classes have been constructed by S. Ivanov (see
\cite[Theorem 41.2]{OlS_book}), and examples with exactly one
nontrivial conjugacy class have recently been announced by D. Osin
\cite{Osin}. But it is an open question to find infinite finitely
presented groups with finitely many conjugacy classes.
\end{rem}

\subsection{Hereditary properties}

\begin{lem}
If $N$ is normal in $G$, then the following statements are
equivalent: (i) $N$ is almost maximal, (ii) $N$ has finite
index.\label{caracfini}\end{lem} \bpr It suffices to show
(i)$\Rightarrow$(ii), which is equivalent to the statement: every
infinite group has infinitely many subgroups. If $G/N$ is torsion,
then it is the union of its finite subgroups, so they are infinitely
many. Otherwise, $G/N$ contains an infinite cyclic subgroup, which
contains infinitely many subgroups.\epr

\begin{lem} Suppose $H_1,H_2$ are subgroups of $G$, $H_1\subset
H_2$. Suppose that $H_1$ has finite biindex (resp. is almost
maximal) in $G$. Then $H_2$ has finite biindex (resp. is almost
maximal) in $G$ and $H_1$ has finite biindex (resp. is almost
maximal) in $H_2$.\label{dicho}\end{lem} \bpr The statement for with
``almost maximal" is trivial.

Suppose that $H_1$ has finite biindex in $G$. Trivially, so has
$H_2$. Write $G=\bigcup_{i\in I}H_1g_iH_1$, with $I$ finite, and set
$J=\{i\in I, g_i\in H_2\}$. Then $H_2=\bigcup_{i\in J}H_1g_iH_1$, so
that $H_1$ has finite biindex in $G$.

\begin{lem}Suppose $H_1,H_2$ are subgroups of $G$, $H_1$ is contained as a subgroup
of finite index in $H_2$, and $H_2$ has finite biindex in $G$. Then
$H_1$ has finite biindex in $G$.\label{fibas}
\end{lem}
\bpr Write $G=\bigcup_i H_2g_iH_2$, $H_2=\bigcup_j h_jH_1=\bigcup_k
H_1h'_k$, all these unions being finite. Then $G=\bigcup_{i,j,k} H_1
h'_kg_ih_jH_1$.\epr

\begin{rem} The converse of Lemma \ref{dicho} is false in both cases.

In the case of finite biindex, consider $G=\textnormal{SL}(2,K)$,
where $K$ is an algebraically closed field, $T$ is the subgroup of
upper triangular matrices in $G$, and $D$ denotes the diagonal
matrices in $G$. Then $G$ is two-transitive on $G/T\simeq P^1(K)$,
and $T$ is two-transitive on $T/D\simeq K$, the affine line. But, by
a dimension argument, the action of $D$ on $G/D$ cannot have a
finite number of orbits. On the other hand, I do not know any
counterexample with $G$ finitely generated.

For a counterexample with almost maximal subgroups, which also shows
that the analogue of Lemma \ref{fibas} is false with almost maximal
subgroups, take an infinite group $G$ with a finite maximal subgroup
$H$. Such groups are constructed in \cite{OlS_book}. So $H$ is
almost maximal in $G$ and $\{1\}$ has finite index in $H$. But, by
Lemma \ref{caracfini} $\{1\}$ is not almost maximal in
$G$.\label{antidicho}\end{rem}

\medskip

\begin{rem}Here is a trivial consequence of Lemma \ref{dicho}.
Let $G_1$ be a group, $G_2$ is a finite index subgroup of $G_1$, and
$H$ a subgroup of $G_2$. Then, if $H$ has finite biindex (resp. is
almost maximal) in $G_2$, it has also finite biindex (resp. is
almost maximal) in $G_1$. The point is that I do not know, in both
cases, if the converse is true.\end{rem}

\begin{lem}
Suppose that, for $i=1,2$, $H_i$ has finite biindex (resp. is almost
maximal) in $G_i$. Then $H_1\times H_2$ has finite biindex (resp. is
almost maximal) in $G_1\times G_2$.
\end{lem}
\bpr This is obvious with finite biindex. Suppose that, for $i=1,2$,
$H_i$ is almost maximal in $G_i$. If there are infinitely many
subgroups containing $H_1\times H_2$, infinitely many have the same
intersection $K_i$ and projection $P_i$ on $G_i$ for $i=1,2$. Note
that $K_i$ is normal in $P_i$. Since, as a consequence of Lemma
\ref{dicho}, $K_i$ is almost maximal in $P_i$, this implies, by
Lemma \ref{caracfini}, that $P_i/K_i$ is finite for $i=1,2$. Since
only finitely many subgroups can exist between $K_1\times K_2$ and
$P_1\times P_2$, we have a contradiction.\epr

\begin{lem}
If $H$ has finite biindex (resp. is almost maximal) in $G$ and $N$
is a normal subgroup of $G$, then $H/(H\cap N)$ has finite biindex
(resp. is almost maximal) in $G/N$.\label{quotient}
\end{lem}
\bpr For the case of finite biindex, pass the expression
$G=\bigcup_{i=1}^n Hg_iH$ to the quotient. The statement for almost
maximal subgroups is trivial.\epr

\begin{prop}
Properties (BF) and (LBF) are inherited from finite index subgroups.
\end{prop}
\bpr Let $G$ be a group and $N$ a finite index subgroup. Suppose
that $N$ is (L)BF. Let $H$ be a (finitely generated) almost maximal
subgroup in $G$. By Lemma \ref{fibas}, $H\cap N$ has finite biindex
in $G$ (and is also finitely generated), so has finite biindex in
$N$ by Lemma \ref{dicho}. Since $N$ has Property (L)BF, $H\cap N$
has finite index in $N$, so that $H$ has finite index in $G$.\epr

\medskip

I do not know if Properties (BF) and (LBF) are inherited by finite
index subgroups. This motivates the following question.

\begin{que}
Let $G$ be a group, $N$ a subgroup of finite index, and $H$ a
subgroup of $N$. If $H$ has finite biindex in $N$, must it have
finite biindex in $G$? \label{Q2}
\end{que}

If Question \ref{Q2} has a positive answer, then Properties (BF) and
(LBF) are inherited by finite index subgroups.

\subsection{Faithful almost 2-transitive pairs}

We could define weaker analogs of Properties (PF) through (LBF), say
$(fPF)$, etc., by only considering subgroups $H$ such that the
action of $G$ on $G/H$ is faithful.

Not much is known about these properties for infinite groups. Dixon
\cite{Dixon} has shown that, in a suitable sense, a generic subgroup
on $n\ge 2$ generators of the symmetric group $\text{Sym}(\N)$ is
free and 2-transitive, showing that free groups also have {\em
faithful} 2-transitive actions on infinite sets and therefore do not
have Property (fBF). The non-existence of a faithful primitive
action, which is, for infinite groups, a priori slightly stronger
than Property (fMF), is widely investigated in \cite{GG}.

Examples \ref{ThompF/FxF}, \ref{ThompF/T}, and \ref{Johnson} provide
essentially the only examples of finitely presented groups which I
know not to have Property (fLBF). In the finitely generated case, we
also have the groups $G\times G$ when the infinite group $G$ has
finitely many conjugacy classes, and has trivial center (this latter
assumption is always satisfied modulo a finite normal subgroup).
Note that, in all these examples, $G$ has very few normal subgroups:
for $F$ and the Houghton groups, $G$ has simple derived subgroup;
$T$ is itself simple, and the groups with finitely many conjugacy
classes have finitely many normal subgroups. We therefore ask:

\begin{que}
Does there exist a residually finite group that acts almost
2-transitively and faithfully on an infinite set, with finitely
generated stabilizers?
\end{que}

The answer is yes when ``almost 2-transitively" is replaced by
``primitively", as shows the example, pointed out in \cite{GG}, of
the action of $\textnormal{PSL}_2(\Z[1/p])$ on
$\textnormal{PSL}_2(\Z[1/p])/\textnormal{PSL}_2(\Z)$.

\appendix

\section{Length of words in wreath products}

We consider the wreath product $A=W\wr_{G/H}G$, where $W$ and $G$
are finitely generated. We write $W$ additively although it is not
necessarily abelian. We write the elements of $A$: $(f,c)$, where
$f\in W^{(G/H)}$ and $c\in G$; we denote by $x_0$ the base-point of
$G/H$. If $w\in W$ and $c\in G$, by the abusive notation $(w,c)$, we
mean $(\delta_{x_0}(w),c)$, where $\delta_{x_0}$ is the natural
inclusion of $W$ into the $x_0$-component of $W^{(G/H)}$. The
product in $A$ is given by $(f_1,c_1)(f_2,c_2)=(f_1+c_1f_2,c_1c_2)$.
Fix a symmetric finite generating subset $S$ of $G$. We call a path
of length $n$ in $G$ a sequence $(g_0,\dots,g_n)$ such that $g_0=1$
and $g_{i+1}g_i^{-1}\in S$ for all $i=0,\dots,n-1$. For any finite
subset $F$ of $X$ and $c\in G$, let $K(F,c)$ be the minimal length
of a path $(g_0,\dots,g_n)$ in $G$ such that $g_n=c$ and
$F\subset\{g_0x_0,\dots,g_nx_0\}$. On the other hand, fix a finite
symmetric generating subset $T$ of $W$, and denote by $|\cdot|$ the
corresponding word length. For $f\in W^{(X)}$, set $|f|=\sum_{x\in
X}|f(x)|$. Fix, as generating subset of $A$, the union of $(t,1)$
for $t\in T$ and $(0,s)$ for $s\in S$. Denote again by $|\cdot|$ the
word length in $A$. The following lemma is obtained in \cite{parry}
in the case of standard wreath products.

\begin{lem}
For $f\in W^{(G/H)}$ and $c\in G$, we have
$|(f,c)|=K(\textnormal{supp}(f),c)+|f|$.
\end{lem}
\bpr Set $n=K(\textnormal{supp}(f),c)$. Let $1=g_0,g_1,\dots,g_n$ be
a path of length $n$ such that $g_n=c$ and, whenever
$x\in\text{Supp}(f)$, $x=g_ix_0$ for some $i$. For all $i$, set
$g_{i+1}=s_ig_i$ with $s_i\in S$, and $x_i=g_ix_0$. Then
$$(f,c)=(f(x_0),s_1)(f(x_1),s_2)\dots (f(x_{n-1}),s_n)(f(x_n),1)$$
$$=(f(x_0),1)(0,s_1)(f(x_1),1)(0,s_2)\dots
(f(x_{n-1}),1)(0,s_n).$$

Thus, $(f,c)$ can be expressed as a product of
$K(\textnormal{supp}(f),c)+|f|$ generators. Accordingly, for all
$(f,c)$, $|(f,c)|\le K(\textnormal{supp}(f),c)+|f|$.

Conversely suppose that $(f,c)$ can be expressed as a product of a
minimal number $n$ of generators. Putting generators of $W_{x_0}$
together, we get
$$(f,c)=(w_1,1)(0,s_1)(w_2,1)(0,s_2)\dots
(w_m,1)(0,s_m)(w_{m+1},1)$$
$$=(w_1,s_1)\dots (w_m,s_m)(w_{m+1},1)$$
with each $w_i\in W$, $s_i\in S$. Set $g_i=\prod_{j=1}^is_i$.

Then $(g_0,g_1,\dots,g_m)$ is a path joining $1$ to $g_m=c$.
Besides, $f=w_1+g_1w_2+\dots+g_mw_{m+1}$, so that $\text{Supp}(f)$
is contained in $\{g_0x_0,g_1x_0,\dots,g_{m-1}x_0,g_mx_0\}$.
Accordingly, $K(\text{Supp}(f),c)\le m$, and $n=\sum_i|w_i|+m\ge
|f|+K(\text{Supp}(f),c)$.\epr

\medskip

This immediately implies the following result, first observed by A.
Erschler-\allowbreak Dyubina \cite{Dyubina} in the case of standard
wreath products.

\begin{prop}
Let $G$ be a group, $H$ a subgroup, and $W_1,W_2$ two
bi-Lipschitz-equivalent groups. Then $W_1\wr_{G/H} G$ and
$W_2\wr_{G/H} G$ are bi-Lipschitz-equivalent.\label{Erschler}
\end{prop}

\bigskip

\noindent\textbf{Acknowledgments.} I am grateful to Laurent
Bartholdi and Luc Guyot for valuable discussions and encouragement.
I thank Victor Guba and Mark Sapir for discussions about Thompson's
groups. I also thank Pierre de la Harpe and Romain Tessera for
useful corrections.

%\bigskip
%\footnotesize
%\noindent Yves de Cornulier\\
%EPF Lausanne, IGAT\\
%E-mail: \url{decornul@clipper.ens.fr}\\
% ----------------------------------------------------------------
\end{document}